\documentclass[12pt]{article}
\usepackage{pgf,tikz}
\usepackage{calc,graphicx, complexity}
\usepackage{graphics}
\everymath{\displaystyle}
\usepackage{graphicx}
\usepackage{color}
\usepackage{amssymb}
\usepackage{amsmath}
\newtheorem{prethm}{{\bf Theorem}}

\newenvironment{thm}{\begin{prethm}{\hspace{-0.5
				em}{\bf .}}}{\end{prethm}}
\newtheorem{prelemma}{{\bf Lemma}}

\newenvironment{lemma}{\begin{prelemma}{\hspace{-0.5
				em}{\bf .}}}{\end{prelemma}}
\newtheorem{preex}{{\bf Example}}

\newtheorem{preprop}{{\bf Proposition}}

\newenvironment{prop}{\begin{preprop}{\hspace{-0.5em}{\bf .}}}{\end{preprop}}
\newtheorem{precor}{{\bf Corollary}}

\newenvironment{cor}{\begin{precor}{\hspace{-0.5
				em}{\bf .}}}{\end{precor}}
\newtheorem{preremark}{{\bf Remark}}

\newtheorem{preprob}{{\bf Problem}}

\newtheorem{predefin}{{\bf Definition}}

\newtheorem{preconj}{{\bf Conjecture}}

\newtheorem{preprobb}{{\bf Problem}}

\newtheorem{prelem}{{\bf Theorem}}

\newtheorem{precla}{{\bf Claim}}

\newenvironment{proof}{{\bf Proof.}\rm }{\hfill{$\Box$}}

\newtheorem{presolution}{{\bf Solution.}}

\def\newpic#1{}
\def\qed{\ifhmode\unskip\nobreak\fi\quad\ifmmode\Box\else$\Box$\fi}

\title{\vspace{-1cm}\Large\bf\noindent On z-coloring and ${\rm b}^{\ast}$-coloring of graphs as improved variants of the b-coloring}
\author{\large\bf Manouchehr Zaker\footnote{mzaker@iasbs.ac.ir}
\vspace{5mm}\\
Department of Mathematics,\\
Institute for Advanced Studies in Basic Sciences,\\
Zanjan 45137-66731, Iran\\
}

\date{}

\begin{document}
\maketitle

\begin{abstract}
\noindent Let $G$ be a simple graph and $c$ a proper vertex coloring of $G$. A vertex $u$ is called b-vertex in $(G,c)$ if all colors except $c(u)$ appear in the neighborhood of $u$. By a ${\rm b}^{\ast}$-coloring of $G$ using colors $\{1, \ldots, k\}$ we define a proper vertex coloring $c$ such that there is a b-vertex $u$ (called nice vertex) such that for each $j\in \{1, \ldots, k\}$ with $j\not=c(u)$, $u$ is adjacent to a b-vertex of color $j$. The ${\rm b}^{\ast}$-chromatic number of $G$ (denoted by ${\rm b}^{\ast}(G)$) is the largest integer $k$ such that $G$ has a ${\rm b}^{\ast}$-coloring using $k$ colors. Every graph $G$ admits a ${\rm b}^{\ast}$-coloring which is an improvement over the famous b-coloring. A z-coloring of $G$ is a coloring $c$ using colors $\{1, 2, \ldots, k\}$ containing a nice vertex of color $k$ such that for each two colors $i<j$, each vertex of color $j$ has a neighbor of color $i$ in the graph (i.e. $c$ is obtained from a greedy coloring of $G$). We prove that ${\rm b}^{\ast}(G)$ cannot be approximated within any constant factor unless $\P=\NP$. We obtain results for ${\rm b}^{\ast}$-coloring and z-coloring of block graphs, cacti, $P_4$-sparse graphs and graphs with girth greater than $4$. We prove that z-coloring and ${\rm b}^{\ast}$-coloring have a locality property. A linear 0-1 programming model is also presented for z-coloring of graphs. The positive results suggest that researches can be focused on ${\rm b}^{\ast}$-coloring (or z-coloring) instead of b-coloring of graphs.
\end{abstract}



\section{Introduction}

\noindent All graphs in this paper are undirected without any loops and multiple edges. In a graph $G$, $\Delta(G)$ denotes the maximum degree of $G$. Let $v$ be a vertex in $G$, $N(v)$ is the set of neighbors of $v$ in $G$. Define the degree of $v$ as $d(v)=|N(v)|$. Also for a subset $S$ of vertices in $G$, by $G[S]$ we mean the subgraph of $G$ induced by the elements of $S$. Complete and path graphs on $n$ vertices are denoted by $K_n$ and $P_n$, respectively. The maximum number of mutually adjacent vertices in $G$ is denoted by $\omega(G)$. The join $G\vee H$ of two vertex disjoint graphs $G$ and $H$ is a graph on $V(G)\cup V(H)$ in which each vertex in $G$ is adjacent to all vertices in $H$. A proper vertex coloring $c$ of a graph $G$ is an assignment of colors $1, 2, \ldots$ to the vertices of $G$ such that no two adjacent vertices receive same colors. By a color class we mean a subset of vertices having a same color. The chromatic number $\chi(G)$ of $G$, is the smallest number of colors used in a proper coloring of $G$. We refer to \cite{BM} for the concepts not defined here. In a proper vertex coloring $c$ of $G$, a vertex $u$ is called b-vertex if $u$ has a neighbor of color $j$ for each color $j\not= c(u)$. A coloring $c$ is b-coloring if for each color $j$, there exists a b-vertex of color $j$. The maximum number of colors in a b-coloring of $G$ is called b-chromatic number and denoted by b$(G)$ (also by $\chi_b(G)$). Clearly, ${\rm b}(G)\leq \Delta(G)+1$. A coloring consisting of color classes say $C_1, \ldots, C_k$ is Grundy-coloring of $G$ if for each $i < j$ each vertex in $C_j$ has a neighbor in $C_i$. The Grundy-coloring can be considered as off-line version of the First-Fit coloring \cite{Z1}. The Grundy number, denoted by $\Gamma(G)$ (also by $\chi_{_{\sf FF}}(G)$) is the maximum number of colors used in a Grundy-coloring of $G$. A vertex $u$ is called nice vertex in \cite{Z2} if for each $i, j$ with $j\not= c(u)$, $u$ has a neighbor which is nice vertex of color $j$. A proper coloring $c$ is called z-coloring in \cite{Z2} if $c$ is a Grundy-coloring using say $k$ colors such that $c$ contains a nice vertex of color $k$. Denote by z$(G)$ the maximum number of colors in a z-coloring of $G$. We have z$(G)\leq \min\{\Gamma(G),{\rm b}(G)\}$, also $\Gamma(G)\leq \Delta(G)+1$. However we can obtain a better bound for $\Gamma(G)$ (hence for z$(G)$) as follows. For $u\in V(G)$, define $\Delta(u)={\max} \{d(v):v\in N(u),~d(v)\leq d(u)\}$ and $\Delta_2(G)={\max}_{u\in V(G)} \Delta(u)$. It was proved in \cite{Z12} that $\Gamma(G)\leq \Delta_2(G)+1$. Hence, z$(G)\leq \Delta_2(G)+1$. Note that $\Delta_2(G)\leq \Delta(G)$ and $\Delta(G)-\Delta_2(G)$ may be arbitrarily large. The literature is full of papers concerning the Grundy number, First-Fit coloring and b-coloring of graphs e.g. \cite{EGT,GL,HK,IM,JP,KPT,KTV,Z1,Z2}. To determine $\Gamma(G)$ is $\NP$-completeness for the complement of bipartite graphs \cite{Z1}. Also to determine ${\rm b}(G)$ is $\NP$-complete for complement of bipartite graphs \cite{BSSV} and bipartite graphs \cite{KTV}. z-coloring of graphs was studied in \cite{KZ,SH,Z2}. It was proved in \cite{Z2} that z$(T)$ can be determined in polynomial time for given trees $T$. Also \cite{KZ} proves that to determine z$(G)$ is $\NP$-complete for bipartite graphs $G$. Let $k$ be a fixed integer. It was proved in \cite{SH} that deciding whether ${\rm z}(G)=k+1$ can be solved by a polynomial time algorithm, where $G$ is a $k$-regular graph.

\noindent For vertex coloring of graphs, many fascinating algorithms, heuristics and metaheuristics have been designed \cite{L}. Although the sophisticated heuristics and metaheuristics are usually more optimal than the algorithms with simple structures but it's not possible to obtain theoretical comparative results for their solutions in comparison to optimal solutions. For example we cannot estimate the outputs of Brelaz's DSATUR and Leighton's Recursive Largest First heuristics even for $P_4$-free graphs (see \cite{L} for definitions). But it can be proved that the simple First-Fit coloring is always optimal for such graphs. Hence, a moderate gaol is to explore heuristics capable of theoretical and analytic abilities and as competitive as possible in the experimental contests. The Grundy-coloring and b-coloring are two widely studied color surpassing procedures (see \cite{Z2}) and optimal to some extents. The z-coloring and ${\rm b}^{\ast}$-coloring (to be defined in this paper) are better heuristics satisfying these properties.

\noindent By a ${\rm b}^{\ast}$-coloring of a graph $G$ using colors $\{1, \ldots, k\}$ we define a proper vertex coloring such that there is a b-vertex $u$ (called nice vertex) of color $k$ such that for each $j<k$, $u$ is adjacent to a b-vertex of color $j$. We call it ${\rm b}^{\ast}$-coloring because the subgraph induced on b-vertices contains a star graph $K_{1,k-1}$ as subgraph, where the nice vertex is at the center. In the context of ${\rm b}^{\ast}$-coloring, a nice vertex can also be called ${\rm b}^{\ast}$-vertex (similar to b-vertex). The ${\rm b}^{\ast}$-chromatic number ${\rm b}^{\ast}(G)$ of $G$ is the largest integer $k$ such that $G$ has a ${\rm b}^{\ast}$-coloring using $k$ colors. The following was proved in \cite{Z2}.

\begin{thm}(\cite{Z2})
Let $G$ be a graph on $n$ vertices and $m$ edges and $c$ a proper vertex coloring of $G$ with color classes $C_1, \ldots, C_k$. Then $c$ can be transformed into a z-coloring with at most $k$ colors using local re-colorings beginning from top class $C_k$ down to $C_1$. Moreover, the transformation takes $\mathcal{O}(nm)$ time steps.\label{z}
\end{thm}

\noindent Based on the algorithm in Theorem \ref{z}, a coloring heuristic called IZ was defined in \cite{Z2} by applying z-coloring procedure iteratively. Theorem \ref{z} implies that every graph $G$ admits a ${\rm b}^{\ast}$-coloring and then ${\rm b}^{\ast}(G)$ is well-defined. We have z$(G)\leq {\rm b}^{\ast}(G) \leq {\rm b}(G)$. It is easily seen that ${\rm b}(G)-{\rm b}^{\ast}(G)$ can be arbitrarily large even for tree graphs. For example ${\rm b}(G)$ can be large for caterpillar trees $G$ but for every such graph we have ${\rm b}^{\ast}(G)\leq 3$. Bonomo et al. introduced the concept of {\rm b}-monotonicity \cite{BDMMV}. A graph $G$ is {\rm b}-monotonic if ${\rm b}(H_2)\leq {\rm b}(H_1)$ for every induced subgraph $H_1$ of $G$ and every induced subgraph $H_2$ of $H_1$. z-monotonic and ${\rm b}^{\ast}$-monotonic graphs are defined similarly. Let $G$ be a graph obtained from removing a matching of size $n$ from $K_{n+1,n+1}$. Then $H=K_{n,n}$ is induced subgraph of $G$ and ${\rm b}^{\ast}(H)={\rm z}(H)=n$ but ${\rm b}^{\ast}(G)=2$. Hence, $G$ is neither z-monotonic nor ${\rm b}^{\ast}$-monotonic.

\noindent In the following we prove that in the definition of z-coloring, a nice vertex should not necessarily have a largest color. This provides a new description of z-coloring. In fact, if we have a coloring $c$ such that $c$ is Grundy-coloring in $G$ and $c$ contains a b-vertex $u$ such that $u$ is adjacent to a b-vertex of each color other than $c(u)$, then $c$ can be transformed to a z-coloring containing a nice vertex having the largest color (Proposition \ref{new-des}). The new description is used to obtain a 0-1 programming model for z-coloring.
For a graph $G$, denote by z$'(G)$ the maximum number $k$ in a Grundy-coloring $c$ of $G$ using $k$ colors such that $c$ contains a b-vertex $u$ of color $c(u)$ such that for each $j\not= c(u)$ with $1\leq j\leq k$, $u$ is adjacent to a b-vertex of color $j$.

\begin{prop}\label{new-des}
For any graph $G$, ${\rm z}'(G)={\rm z}(G)$.
\end{prop}

\noindent \begin{proof}
Obviously ${\rm z}'(G)\geq {\rm z}(G)$. To prove the inverse inequality, let z$'(G)=k$ and $c$ consisting of classes $C_1, \ldots, C_k$ be a Grundy-coloring of $G$ using $k$ colors containing a nice vertex $u$. Represent the color classes from down to up, i.e. $C_1$ is the lowest class and $C_k$ the most top one. If $u\in C_k$ then ${\rm z}(G)\geq {\rm z}'(G)$. Let $u\in C_j$, $j<k$. Consider a new proper coloring $c'$ of $G$ as follows. Re-color all vertices in $C_j$ by $k$. Denote this class in $c'$ by $C'_k$ which is now the most top class in $c'$. For each $i$ with $j+1\leq i \leq k$, re-color all vertices in $C_i$ by $i-1$. Clearly, in the coloring $c'$, the color classes lower than $C'_k$ induce a Grundy-coloring. Each b-vertex in $c$ remains b-vertex in $c'$. We Grundyfy the class $C'_k$. In other words, transfer each vertex $w$ in $C'_k$ to the lowest class in which $w$ has not any neighbor. Since $u\in C'_k$ and $u$ has a b-vertex neighbor in $c'$ then the color of $u$ remains $k$ after Grundyfing the class $C'_k$. It means that if $C''_k$ denotes the class $C'_k$ after Grundying then $C''_k\not= \emptyset$. It follows that $c'$ is a Grundy-coloring using $k$ colors containing a nice vertex of color $k$. Then by the definition z$(G)\geq k$. Then z$'(G)$$=$z$(G)$.
\end{proof}

\noindent For a graph $G$ of degree sequence $d_1\geq \cdots \geq d_n$, define $m(G)=\max \{i: d_i\geq i-1\}$. It is known that ${\rm b}(G)\leq m(G)$. The invariant can be generalized for ${\rm b}^{\ast}$-coloring. Define $m^{\ast}(G)$ as the maximum $k$ such that there exists $u\in V(G)$ and $u_1, \ldots, u_k \in N(u)$ such that for each $i$, $d(u_i)\geq k$. It is easily seen that ${\rm b}^{\ast}(G) \leq m^{\ast}(G)+1$. For this purpose, let $c$ be a ${\rm b}^{\ast}$-coloring of $G$ using say $k={\rm b}^{\ast}(G)$ colors. Let $u$ be a nice vertex of color $k$ in $c$. Vertex $u$ has at least $k-1$ neighbors each of degree at least $k-1$. It follows by the definition that $m^{\ast}(G)\geq k-1$. We have also
$${\rm z}(G) \leq {\rm b}^{\ast}(G) \leq m^{\ast}(G)+1 \leq \Delta_2(G)+1 \leq \Delta(G)+1.$$
\noindent Note that $m^{\ast}(G)-{\rm b}^{\ast}(G)$ can be arbitrarily large. Let $H_n$ be a graph obtained by removing a matching of size $n-1$ from the complete bipartite graph $K_{n,n}$. We have ${\rm b}^{\ast}(H_n)\leq {\rm b}(H_n)\leq 2$ but $m^{\ast}(H_n)=n-1$.

\noindent Proposition \ref{m-algo} determines $m^{\ast}(G)$ by an efficient algorithm. We need some knowledge about a fast sorting algorithm. Given a set of positive integers $B$ of cardinality $k$ such that $b\leq k$, for each $b\in B$. Using the counting sort (see Page 194 in \cite{CLRS}) we can sort the elements of $B$ non-increasingly with time complexity $\mathcal{O}(k)$. As a bypass result, we can determine using $\mathcal{O}(k)$ operations the maximum $p$ such that there exists $b_1, \ldots, b_p$ satisfying $b_i\geq p$, for each $i=1, \ldots, p$. For this purpose sort $B$ non-increasingly such as $b_1 \geq \cdots \geq b_k$ and obtain with $k$ comparisons the largest index $p(B)$ such that $b_i\geq i$, for $i=1, \ldots, p$ and $b_{p+1}\leq p+1$. This $p(B)$ is the solution.

\begin{prop}\label{m-algo}
There exists an $\mathcal{O}(n\Delta)$ time algorithm which determines $m^{\ast}(G)$ for any graph $G$ on $n$ vertices with maximum degree $\Delta$.
\end{prop}

\noindent \begin{proof}
Let $G$ on $n$ vertices be presented by its list of adjacency. For each $u\in V(G)$, define $A_u=\{d(w):w\in N(u)\}$. By the previous paragraph we obtain $p(A_u)$ by consuming $\mathcal{O}(\Delta(G))$ time steps. It is easily seen that $m^{\ast}(G)={\max}_{u\in G} p(A_u)$. It follows that $m^{\ast}(G)$ is determined with time complexity $\mathcal{O}(n\Delta(G))$.
\end{proof}

\noindent Also we have the following proposition which does not hold for b-chromatic number.

\begin{prop}\label{disj-union}
Let $G$ be a vertex disjoint union of $G_1, \ldots, G_p$. Then ${\rm z}(G)=\max_{i=1}^p {\rm z}(G_i)$ and ${\rm b}^{\ast}(G)=\max_{i=1}^p {\rm b}^{\ast}(G_i)$.
\end{prop}

\noindent {\bf The rest of the paper is organized as follows.} We prove in Section 2 that ${\rm b}(G)={\rm b}^{\ast}(G \vee K_1)-1$ (Proposition \ref{eq}). Proposition \ref{hard} asserts that ${\rm b}^{\ast}(G)$ cannot be approximated within any constant ratio unless $\P=\NP$. Also to determine ${\rm b}^{\ast}(G)$ is $\NP$-hard for co-bipartite graphs. Then we prove that block and cactus graphs are {\rm z}-monotonic and ${\rm b}^{\ast}$-monotonic (Propositions \ref{block}). It follows that
${\rm b}^{\ast}(G)=m^{\ast}(G)+1$ for block graphs $G$ (Proposition \ref{block=}). A similar result is obtained for graphs of girth five in Theorem \ref{girth}. Then Proposition \ref{zsparse} proves that ${\rm z}(G)=\omega(G)$ for $P_4$-sparse graphs $G$. In Section 3, we prove a locality result for z-chromatic and ${\rm b}^{\ast}$-chromatic numbers in Proposition \ref{zlocal}. A linear 0-1 programming model for z$(G)$ is presented and proved in Proposition \ref{prog}. The paper ends with introducing some open research areas.

\section{Results on ${\rm b}^{\ast}(G)$ and z$(G)$}

\noindent Let $G$ be a graph. By $G \vee K_1$ we mean a graph obtained by adding an extra vertex to $G$ and connecting it to each vertex in $G$.

\begin{prop}
For any graph $G$, ${\rm b}(G)={\rm b}^{\ast}(G \vee K_1)-1$.\label{eq}
\end{prop}

\noindent \begin{proof}
Write $H=G \vee K_1$ and denote by $w$ the vertex of $H$ joined to $V(G)$. Let $c$ be a {\rm b}-coloring of $G$ using $k={\rm b}(G)$ colors and with $u_1, \ldots, u_k$ as its b-vertices. Define a coloring $c'$ for $H$ as follows. Set $c'(v)=c(v)$, for each $v\in V(G)$. Also $c'(w)=k+1$. Clearly, $c'$ is a ${\rm b}^{\ast}$-coloring of $H$ in which $w$ is a nice vertex since each $u_i$ is b-vertex in $c'$. It follows that ${\rm b}^{\ast}(H)\geq {\rm b}(G)+1$. Now let $c''$ be a ${\rm b}^{\ast}$-coloring of $H$ using $t={\rm b}^{\ast}(H)$ colors. The color class in $c''$ containing $w$ consists only $w$. W.l.o.g we may assume that $c''(w)=t$. There are b-vertices say $v_1, \ldots, v_{t-1}$ other than $w$ in $c''$. These vertices are b-vertex in $G$ with the coloring of $G$ obtained by restriction of $c''$ on $V(G)$. It follows that ${\rm b}(G)\geq t-1$. This completes the proof.
\end{proof}

\noindent There are many hardness and inapproximability results for the b-coloring problem \cite{BSSV,CVV,GK,KTV}. Is was proved in \cite{BSSV} that to determine b$(G)$ is $\NP$-complete for complement of bipartite graphs (shortly co-bipartite graphs). It was proved in \cite{GK} that, for all $\epsilon >0$, it is $\NP$-hard to approximate the b-coloring problem for graphs with $n$ vertices within a factor $n^{(1/4) -\epsilon}$.

\begin{prop}
To determine the ${\rm b}^{\ast}$-chromatic number of co-bipartite graphs is $\NP$-complete. Also for all $\epsilon >0$, it is $\NP$-hard to approximate the ${\rm b}^{\ast}$-coloring problem for graphs with $n$ vertices within a factor $n^{(1/4) - \epsilon}$. In particular, no polynomial-time approximation algorithm within any constant ratio exists for ${\rm b}^{\ast}(G)$ unless $\P=\NP$.\label{hard}
\end{prop}

\noindent \begin{proof}
We transform an instance $(G,k)$ of the b-coloring problem into $(H,k+1)$ of the ${\rm b}^{\ast}$-coloring problem, where
$k$ is an integer and $H=G \vee K_1$. By Proposition \ref{eq}, ${\rm b}^{\ast}(G \vee K_1)-{\rm b}(G)=1$. Note that if $G$ is co-bipartite graph then $H$ is co-bipartite too. It follows by the previous paragraph (i.e. a result of \cite{BSSV}) that to determine ${\rm b}^{\ast}(G)$ is $\NP$-complete for co-bipartite graphs.

\noindent The equality ${\rm b}^{\ast}(G \vee K_1)-{\rm b}(G)=1$ means that the transformation is gap-preserving. Hence,
every approximation hardness result for the b-coloring problem such as the one within the factor $n^{(1/4) -\epsilon}$ mentioned in the previous paragraph and proved in \cite{GK} holds for ${\rm b}^{\ast}$-coloring problem too.
\end{proof}

\noindent b-colorings of block graphs is the subject of many papers e.g. \cite{HSS}. We don't know whether or not ${\rm b}(G)$ is polynomially computable for block graphs. But ${\rm b}^{\ast}(G)$ is easily determined by Proposition \ref{block} for block graphs. We need a helping lemma.

\begin{lemma}
Let $G$ be a block graph and $B$ a block of $G$ on $t$ vertices. Assign colors $1, \ldots, t$ arbitrarily to the vertices of $B$. Then the pre-coloring can be extended to a Grundy-coloring of $G$ using $\omega(G)$ colors.\label{exten}
\end{lemma}

\noindent \begin{proof}
The proof is by induction on the number of vertices. Let $B'$ be an arbitrary block in $G$ intersecting $B$ and let $u\in V(B)\cap V(B')$. Let $j$ be the color of $u$ in the pre-coloring of $B$.

\noindent {\bf Case 1. $j\leq |B'|$}

\noindent In this case, assign colors $1, 2, \ldots, |B'|$ to the vertices of $B'$ such that $u$ receives color $j$.
Let $H$ be a connected component of $G\setminus (B\setminus u)$ which contains $B'$. Apply the induction hypothesis for $H$, $B'$ and the coloring of vertices in $B'$ by $1, 2, \ldots, |B'|$ and obtain a Grundy-coloring of $H$ with $\omega(H)\leq \omega(G)$ colors consisting with the colors of $B'$.

\noindent {\bf Case 2. $j> |B'|$}

\noindent In this case, let $H'$ be a connected component of $G\setminus u$ containing $B'\setminus u$.
We have not any color constrains on the vertices of $B'\setminus u$. A minimum coloring of $H'$ using $\omega(H')$ colors can easily be transformed to a Grundy-coloring, as desired.

\noindent We repeat the above argument for other blocks in $G$ intersecting $B$ and obtain a desired Grundy-coloring of $G$ with $\omega(G)$ colors extending the colors in $B$.
\end{proof}

\begin{prop}\label{block}
Let $G$ be a block graph. Then $G$ is {\rm z}-monotonic and ${\rm b}^{\ast}$-monotonic.
\end{prop}

\noindent \begin{proof}
Let $G$ be a block graph. We prove for z-number. The proof for ${\rm b}^{\ast}$-number is similar.
Let $u$ be an arbitrary vertex in $G$. Write for simplicity $G'=G\setminus u$. We prove that ${\rm z}(G')\leq {\rm z}(G)$.

\noindent {\bf Case 1.} $u$ is not cut-vertex and belongs to a block of $G$ say $B$ on $t$ vertices.

\noindent Let $c$ be a z-coloring of $G'$ with $k={\rm z}(G')$ colors. In this case ${\rm z}(G')\geq \omega(G')\geq t-1$. If ${\rm z}(G')=t-1$ then ${\rm z}(G)\geq \omega(G) = t > {\rm z}(G')$. If ${\rm z}(G')\geq t$ and $c$ uses the colors $1, \ldots, k$, for some $k\geq t$, then a Grundy extension of $c$ makes $u$ use a color from $\{1, \ldots, k\}$, since $u$ has $t-1<k$ neighbors in $G'$. Hence $c$ is extended to a z-coloring of $G$ using $k$ colors, i.e. ${\rm z}(G)\geq k$.

\noindent {\bf Case 2.} $u$ is cut-vertex in $G$.

\noindent In this case, let $G_1, \ldots, G_p$ be the vertex disjoint connected components in $G'$. By Proposition \ref{disj-union}, there exists $i$ such that ${\rm z}(G_i)={\rm z}(G')$. Assume w.l.o.g that ${\rm z}(G_1)={\rm z}(G')=k$. Then there exists a z-coloring $c$ of $G'$ such that a nice vertex and its $k-1$ b-vertex neighbors belong to $V(G_1)$. Let $B_1, \ldots, B_p$ be the blocks containing $u$. Hence, $B'_i= B_i\setminus u$ is a block in $G_i$, $i=1, \ldots, p$.

\noindent {\bf Subcase 2.1.} ${\rm z}(G_1)\leq \omega(G)$.

\noindent In this case we have ${\rm z}(G')={\rm z}(G_1)\leq \omega(G)\leq {\rm z}(G)$, as desired.

\noindent {\bf Subcase 2.2.} ${\rm z}(G_1)\geq \omega(G)+1$.

\noindent In this case let $c'$ be the restriction of $c$ on $V(G_1)$. By our notation, the neighbors of $u$ in $G_1$ forms the complete subgraph $B'_1$. The number of colors used in $c'$, i.e. $k$ is greater than $|B'_1|$. We choose a smallest available color from $\{1, \ldots, k\}$ and assign it to $u$. Denote this color by $c^{\ast}(u)$. In the rest of proof, for each $i=2, \ldots, p$, we obtain a Grundy-coloring of $G_i\cup \{u\}$, denoted by $c_i$, using at most $k$ colors in which the color of $u$ is $c^{\ast}(u)$, i.e. $c_i(u)=c^{\ast}(u)$. The contamination of $c', c_2, \ldots, c_p$ is a z-coloring of $G$ with $k$ colors and $c^{\ast}(u)$ as the color of $u$. This will complete the proof. For this purpose, fix an arbitrary $i$, $2\leq i \leq p$. The neighbors of $u$ in $G_i$ is the complete subgraph $B_i\setminus u$. Recall that $k={\rm z}(G_1)\geq \omega(G)+1$. There are two possibilities for $c^{\ast}(u)$:

\noindent If $c^{\ast}(u)\leq |B_i|-1$, then assign colors $1, 2, \ldots, c^{\ast}(u)-1$ arbitrarily to the some vertices in $B'_i$, assign $c^{\ast}(u)$ to $u$ in $B_i$. Now, apply Lemma \ref{exten} for the block graph $G_i\cup \{u\}$ and pre-coloring $V(B_i)$. Obtain a Grundy-coloring of $G_i\cup \{u\}$ with at most $\omega(G)$ colors which is consistent with the colors in $B_i$. Name the coloring $c_i$. This complete the proof for this possibility.

\noindent But if $c^{\ast}(u)\geq |B_i|$, we only take a Grundy-coloring of $G_i$ using at most $\omega(G)$ colors. By adding $u$ and its color $c^{\ast}(u)$ to the coloring, we obtain a desired coloring $c_i$.
\end{proof}

\begin{prop}
For any block graph $G$, ${\rm b}^{\ast}(G)=m^{\ast}(G)+1$. In particular if $T$ is a tree then ${\rm b}^{\ast}(T)=m^{\ast}(T)+1$.\label{block=}
\end{prop}

\noindent \begin{proof}
Let $m=m^{\ast}(G)$. There exists a vertex $u$ which has $m$ neighbors $w_1, \ldots, w_m$ of degree at least $m$. Since $G$ is block graph then for any neighbors $x$ and $y$ of $u$, either $N(x)\cap N(y)=\{u\}$ or $(N(x)\cap N(y))\cup \{u,x,y\}$ is contained in a same block. Let $n_i$ neighbors of $u$ belong to a block $B_i$, $i=1, 2, \ldots, p$, for some $p\geq 1$ such that $n_1+n_2+\cdots+n_p=m$. We define a partial ${\rm b}^{\ast}$-coloring of $G$ using $m+1$ colors. Set $N_i=N(u)\cap B_i$, for $i=1, \ldots, p$. Assign colors $1, 2, \ldots, n_1$ to vertices in $N_1$. Assign colors $n_1+1, \ldots, |B_1|-1$ to the rest of vertices in $B_1\setminus u$. Each vertex say $w_j$ in $N_1$ has $|B_1|-n_1-1$ neighbors such that no one is adjacent to $w_i$ for each $i\not= j$. We assign suitable colors to these neighbors in order to make $w_j$ a b-vertex. Then assign colors $n_1+1, \ldots, n_1+n_2$ to vertices in $N_2$. Then assign distinct colors from $\{1, \ldots, m\}\setminus \{n_1+1, \ldots, n_1+n_2\}$ to the rest of vertices in $B_2 \setminus u$. Then similar to the vertices in $N_1$ assign suitable colors to each $w_j$ in order to make it b-vertex. By repeating this method for all vertices in $N_3\cup \cdots \cup N_p$, we obtain $m$ b-vertices of distinct colors $1, \ldots m$. Assign finally color $m+1$ to $u$. The resulting pre-${\rm b}^{\ast}$-coloring of $G$ is extended to a ${\rm b}^{\ast}$-coloring of whole $G$ with $m+1$ colors. Therefore  ${\rm b}^{\ast}(G)=m^{\ast}(G)+1$.
\end{proof}

\noindent It was proved in \cite{CSMS} that if $G$ is a connected cactus and $m(G)\geq 7$, then ${\rm b}(G)\geq m(G)-1$. So for these cacti, ${\rm b}(G)\in\{m(G)-1,m(G)\}$. But we have a better result for the ${\rm b}^{\ast}$-chromatic number of cactus graphs.

\begin{prop}\label{cactus}
Let $G$ be a cactus graph. Then $G$ is z-monotonic and ${\rm b}^{\ast}$-monotonic. Also ${\rm b}^{\ast}(G)=m^{\ast}+1$.
\end{prop}

\noindent \begin{proof}
We prove for ${\rm b}^{\ast}$-number. The proof for ${\rm z}(G)$ is similar. We prove that ${\rm b}^{\ast}(G\setminus v) \leq {\rm b}^{\ast}(G)$. If ${\rm b}^{\ast}(G)\leq 2$ then obviously the desired inequality holds. Assume that ${\rm b}^{\ast}(G)\geq 3$. If $v$ is not cut-vertex then $d(v)=2$. Let $c$ be a ${\rm b}^{\ast}$-coloring of $G$ using $k\geq 3$ colors. Obviously a minimum color has not appeared in the neighborhood of $v$. By assigning this color to $v$ we obtain a ${\rm b}^{\ast}$-coloring of $G$ using $k$ colors. Assume now that $v$ is cut-vertex in $G$. Let $H$ be a connected component of $G\setminus v$. The two neighbors $x, y$ of $v$ in $H$ have degree one. Hence, their colors in any Grundy-coloring or ${\rm b}^{\ast}$-coloring $c$ of $H$ are $1$ or $2$. It follows that there exists a minimum color other than $1$ and $2$ to color $v$. Repeat this argument for other components in $G\setminus v$. In this case too we obtain a ${\rm b}^{\ast}$-coloring of $G$ using $k$ colors.

\noindent To prove ${\rm b}^{\ast}(G)\geq m^{\ast}+1$, let $u$ be a cut-vertex with $m=m^{\ast}$ neighbors $u_1, \ldots, u_m$ each of degree at least $m^{\ast}$. Assign color $i$ to $u_i$, $i=1, \ldots, m$. Let $u_j$ be such that no block containing $u_j$ contains another $u_i$. In this case the colors $1, \ldots, j-1, j+1, \ldots, m$ can be assigned to the neighbors of $u_j$ which makes $u_j$ a b-vertex of color $j$. Note that these neighbors of $u_j$ have distance at least three from other $u_i$ vertices.

\noindent Consider the case where some $u_j$ and $u_t$ belong to a same block $B$. Either $u_j$ and $u_t$ are adjacent and their other neighbors have distance at least three from each other or $u_j$ and $u_t$ have a unique common neighbor in $B$ and are not adjacent and their other neighbors have distance at least three from each other. In either case we can assign a color to common neighbor of $u_j$ and $u_t$ (if any) and colors $1, \ldots, j-1, j+1, \ldots, m$ (resp. $1, \ldots, t-1, t+1, \ldots, m$) to the neighbors of $u_j$ (resp. $u_t$). Under this pre-coloring $u_i$ is b-vertex of color $i$ and then $u$ is ${\rm b}^{\ast}$-vertex of color $m+1$. Since the graph is ${\rm b}^{\ast}$-monotonic then ${\rm b}^{\ast}(G)\geq m^{\ast}+1$. This completes the proof.
\end{proof}

\noindent A graph without induced subgraph isomorphic to $P_4$ is called cograph. A graph $G$ is $P_4$-sparse if no subgraph on five vertices induces more than one $P_4$.
A spider is a graph with vertices $c_1, \ldots, c_k, s_1, \ldots, s_k$, $k\geq 2$ such that the vertices $c_i$’s form a
clique, the vertices $s_i$’s form an independent set, each $c_i$ is adjacent only to $s_i$; there may be a vertex
$c_0$ that is adjacent to all $c_i$’s and to no $s_i$ ($c_0$ may or may not be present). Spider graphs can be recognized in linear time \cite{JO}.

%A spider is a graph whose vertex set can be partitioned into $C, S$ and $R$, where $C = \{c_1, \ldots, c_k\}$ is a complete set; $S =\{s_1, \ldots, s_k\}$ ($k \geq 2$) is a stable set; ; $R$ is allowed to be empty and if it is not, then all the vertices in $R$ are adjacent
%to all the vertices in $C$ and non-adjacent to all the vertices in $S$. Clearly, the complement
%of a thin spider is a thick spider, and viceversa. The triple $(C, S, R)$ is called the spider partition, and can be found in linear time \cite{JO}.

\begin{thm}(\cite{JO,HK})
If $G$ is a $P_4$-sparse graph then $G$ or $\bar{G}$ is disconnected, or $G$ or $\bar{G}$ is a spider.
%If $G$ is a non-trivial $P_4$-sparse graph, then either $G$ or $\bar{G}$ is not connected, or $G$ is a spider.
\label{spider}
\end{thm}

\noindent The b-chromatic number of cographs and $P_4$-sparse graphs was studied by Bomono et al. \cite{BDMMV}. Although there are interesting results, but no upper bounds in terms of the clique number $\omega(G)$ for b$(G)$ of $P_4$-free of $P_4$-sparse graphs $G$ have been known in the literature. It was proved in \cite{HK} that if $G$ is $P_4$-free then ${\rm b}(H)=\omega(H)$ for any induced subgraph $H$ of $G$ if and only if $G$ does not contain two certain graphs as induced subgraph. A similar result was obtained in that paper for $P_4$-sparse graphs. But the situation is much better for the ${\rm z}$-chromatic number. In the following we use a result from \cite{KZ} that ${\rm z}(G\vee H)=
{\rm z}(G)+{\rm z}(H)$, where $G \vee H$ is the join graph of $G$ and $H$. 
 
\begin{prop}
Let $G$ be a $P_4$-sparse graph. Then ${\rm z}(G)=\omega(G)$.\label{zsparse}
\end{prop}

\noindent \begin{proof}
First, let $G$ be a spider graph. Denote $C = \{c_1, \ldots, c_k\}$ and $S =\{s_1, \ldots, s_k\}$. Call $G$ thin spider if  $s_i$ is adjacent to $c_j$ if and only if $i = j$. Call $G$ thick spider if $s_i$ is adjacent to $c_j$ if and only if $i\not= j$. Note that complement of a thin spider (with or without $c_0$) is thick spider and vice versa. Let $G$ be a thin spider. Then degree of every vertex in $C$ is $|C|+1$ (resp. $|C|$) if $c_0$ exists (resp. $c_0$ does not exist). Also degree of every vertex in $S$ is one in thin spider graphs. It follows that $m^{\ast}(G)=|C|=\omega(G)-1$, if $c_0$ exists, and $m^{\ast}(G)=|C|-1=\omega(G)-1$, if $c_0$ does not exist. These relations imply that ${\rm z}(G)=\omega(G)$.

\noindent If $G$ is thick spider then degree of every vertex in $C$ is $2|C|-1$ (resp. $2|C|-2$) if $c_0$ exists (resp. $c_0$ does not exist). Also degree of every vertex in $S$ is $|C|-1$. Then $m^{\ast}(G)=\omega(G)-1$. It follows that ${\rm z}(G)=\omega(G)$ for thick spiders $G$.

\noindent Rest of the proof is by induction on $|V(G)|$. Note that complement of a thin spider is thick and vice versa.  By Theorem \ref{spider} it suffices to prove for the cases where $G=G_1 \cup G_2$ or $\bar{G}=H_1\cup H_2$. In the first case, by Proposition \ref{disj-union} and applying the induction hypothesis for $G_1$ and $G_2$, we have
$${\rm z}(G)=\max\{{\rm z}(G_1),{\rm z}(G_2)\}=\max\{\omega(G_1),\omega(G_2)\}=\omega(G).$$
\noindent In the second case $G=H_1 \vee H_2$, where $\vee$ is the join notation. We have ${\rm z}(G)={\rm z}(H_1)+{\rm z}(H_2)$. Applying now the induction hypothesis for $H_1$ and $H_2$ we obtain
$${\rm z}(G)={\rm z}(H_1)+{\rm z}(H_2)=\omega(H_1)+\omega(H_2)=\omega(G).$$
\end{proof}

\noindent In order to state the next result, we need to define an $m\times (m-1)$ array $A_m$ on the entry set $\{1, \ldots, m\}$, for each $m\geq 2$. The rows (resp. columns) are indexed by $1, \ldots, m$ (resp. $1, \ldots, m-1$) from up to down (resp. from left to right). Denote by $A(i,j)$ the entry of $A_m$ in row $i$ and column $j$. The array is so that $A(i,j)=A(j+1,i)$ for each $i<m$ and $j<m$. For each $j$, $1\leq j \leq m-1$, the entries in the diagonal $\{(1,j), (2, j+1), (3, j+2), \ldots\}$ are defined as $A(1,j)=j$ and $A(t+1,j+t)=A(t,j+t-1)+(m-t)$, for each $t=2, 3, \ldots$. Arrays $A_4$ and $A_7$ are depicted in Figure \ref{fig1}. Every two rows in $A_m$ have exactly one common entry and for each entry $j=1, \ldots, m$, there are two unique rows containing $j$. \\

\begin{figure}[h]
	\[
\begin{tabular}{|ccc|}
  \hline
  $\fbox{1}$ & $2$ & $3$ \\[0.5eM]
  \hline
  $1$ & $\fbox{4}$ & 5 \\[0.5eM]
  \hline
  2 & 4 & $\fbox{6}$ \\[0.5eM]
  \hline
  3 & 5 & 6 \\[0.5eM]
  \hline
\end{tabular}
\hspace*{2cm}
\begin{tabular}{|cccccc|}
  \hline
  $\fbox{1}$ & $2$ & $3$ & $4$ & $5$ & $6$ \\[0.64eM]
  \hline
  $1$ & $\fbox{7}$ & 8 & 9 & 10 & 11 \\[0.65eM]
  \hline
  2 & 7 & $\fbox{12}$ & 13 & 14 & 15 \\[0.65eM]
  \hline
  3 & 8 & 12 & $\fbox{16}$ & 17 & 18 \\[0.65eM]
  \hline
  4 & 9 & 13 & 16 & $\fbox{19}$ & 20 \\[0.65eM]
  \hline
  5 & 10 & 14 & 17 & 19 & $\fbox{21}$ \\[0.65eM]
  \hline
  6 & 11 & 15 & 18 & 20 & 21 \\[0.6eM]
  \hline
\end{tabular}
\]
\caption{Arrays $A_4$ (left) and $A_7$ (right)}\label{fig1}
\end{figure}

\noindent The ${\rm b}$-chromatic number of graphs in terms of girth is the research subject of many papers, e.g. \cite{CLS}. The following result concerns ${\rm b}^{\ast}$-coloring of graphs of girth at least five.

\begin{thm}
Let $G$ be a ${\rm b}^{\ast}$-monotonic graph of girth greater than four such that no two 5-cycles intersect in a path of length two. Then ${\rm b}^{\ast}(G)=m^{\ast}(G)+1$.\label{girth}
\end{thm}

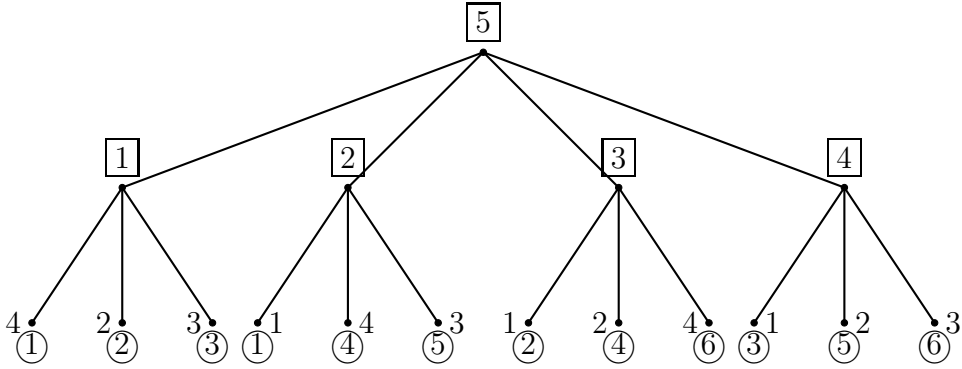
\begin{figure}
\unitlength=0.25mm
\linethickness{0.5pt}
	\begin{center}
		\begin{tikzpicture}[scale=0.60]			
\draw[black, thick] (-8,3)-- (-6,0);
\draw[black, thick] (-8,3)-- (-10,0);	
\draw[black, thick] (-8,3)-- (-8,0);
\draw[black, thick] (-8,3)-- (0,6);

\draw[black, thick] (-3,3)-- (0,6);
\draw[black, thick] (-3,3)-- (-5,0);
\draw[black, thick] (-3,3)-- (-3,0);
\draw[black, thick] (-3,3)-- (-1,0);

\draw[black, thick] (3,3)-- (0,6);
\draw[black, thick] (3,3)-- (1,0);
\draw[black, thick] (3,3)-- (3,0);
\draw[black, thick] (3,3)-- (5,0);

\draw[black, thick] (8,3)-- (0,6);
\draw[black, thick] (8,3)-- (8,0);
\draw[black, thick] (8,3)-- (10,0);
\draw[black, thick] (8,3)-- (6,0);

\filldraw [black] (-10,0) circle(2pt)
node [anchor=north]{$\textcircled{1}$} node [anchor=east]{$4$};
\filldraw [black] (-8,0) circle(2pt)
node [anchor=north]{$\textcircled{2}$} node [anchor=east]{$2$};
\filldraw [black] (-6,0) circle(2pt)
node [anchor=north]{$\textcircled{3}$} node [anchor=east]{$3$};

\filldraw [black] (-5,0) circle(2pt)
node [anchor=north]{$\textcircled{1}$} node [anchor=west]{$1$};
\filldraw [black] (-3,0) circle(2pt)
node [anchor=north]{$\textcircled{4}$} node [anchor=west]{$4$};
\filldraw [black] (-1,0) circle(2pt)
node [anchor=north]{$\textcircled{5}$} node [anchor=west]{$3$};		

\filldraw [black] (1,0) circle(2pt)
node [anchor=north]{$\textcircled{2}$} node [anchor=east]{$1$};		
\filldraw [black] (3,0) circle(2pt)
node [anchor=north]{$\textcircled{4}$} node [anchor=east]{$2$};	
\filldraw [black] (5,0) circle(2pt)
node [anchor=north]{$\textcircled{6}$} node [anchor=east]{$4$};		

\filldraw [black] (6,0) circle(2pt)
node [anchor=north]{$\textcircled{3}$} node [anchor=west]{$1$};		
\filldraw [black] (8,0) circle(2pt)
node [anchor=north]{$\textcircled{5}$} node [anchor=west]{$2$};		
\filldraw [black] (10,0) circle(2pt)
node [anchor=north]{$\textcircled{6}$} node [anchor=west]{$3$};		

\filldraw [black] (0,6) circle(2pt)
node [anchor=south]{$\fbox{5}$};			
\filldraw [black] (-8,3) circle(2pt)
node [anchor=south]{$\fbox{1}$};
\filldraw [black] (-3,3) circle(2pt)
node [anchor=south]{$\fbox{2}$};
\filldraw [black] (3,3) circle(2pt)
node [anchor=south]{$\fbox{3}$};
\filldraw [black] (8,3) circle(2pt)
node [anchor=south]{$\fbox{4}$};
\end{tikzpicture}
\end{center}
\caption{A ${\rm b}^{\ast}$-coloring with five colors obtained from the table in Figure \ref{fig1} (left)}
\label{fig2}
\end{figure}

\noindent \begin{proof}
Let $m=m^{\ast}(G)$. Then there exists $u\in V(G)$ and $m$ neighbors $u_1, \ldots, u_m$ of $u$ such that $d_G(u_i)\geq m$. For each $i=1, \ldots, m$, take a set $S_i$ containing $m-1$ neighbors of $u_i$ other than $u$. Since the girth is at least five then $S_i\cap S_j=\emptyset$, for each $i,j$ with $i\not= j$. Define $H=G[S_1 \cup \cdots \cup S_m]$, the subgraph of $G$ induced by $S_1 \cup \cdots \cup S_m$. The properties of $G$ imply that there exists at most one edge between any two $S_i$ and $S_j$. W.l.o.g. we may assume that there exists exactly one edge between each $S_i$ and $S_j$, $i\not= j$. We have $|V(H)|=m(m-1)$ and $|E(H)|=m(m-1)/2$. Also every vertex of $H$ has degree at most one, otherwise we obtain two $5$-cycles intersecting in a path of length $2$.

\noindent We need a representation of $H$ by an array of integers. Consider a labeling $\ell: V(H) \rightarrow \{1, \ldots, m(m-1)/2\}$ such that for each $x,y\in V(H)$, $xy\in E(H)$ if and only if $\ell(x)=\ell(y)$. Since $|E(H)|=m(m-1)/2$ then there are exactly $m(m-1)/2$ labels in $H$ corresponding to the edges of $H$. We represent $H$ by an array $A_m$ of size $m\times (m-1)$ which displays vertex labels and in fact the edges in $H$. Arrays $A_4$ and $A_7$ are displayed in Figure \ref{fig1}. Let the entries in $i$-th row of $A_m$, denote the vertex labels in $S_i$ (with an arbitrary ordering of vertices in $S_i$). Two vertices in $H$ are adjacent if and only if they have same labels in $A_m$. Let $D$ and $C_j$ ($1\leq j\leq m-1$) be the sets consisting of entries in the main diagonal and $j$-th column of $A_m$, respectively. The elements of $D$ are specified by boxed entries in $A_4$ and $A_7$ in Figure \ref{fig1}. Array $A_m$ provides information which is used to obtain a ${\rm b}^{\ast}$-coloring of $H\cup \{u, u_1, \ldots, u_m\}$ and then $G$ with $m+1$ colors. Each entry in $A_m$ corresponds to a vertex in $H$.

\noindent We show that the entries in $A_m$ are partitioned into $m$ subsets of cardinality $m-1$ such that the elements in each subset are all distinct. The partition sets are $C_1\setminus D, \ldots, C_{m-1}\setminus D$. In fact $C_1$ is the first partition set and $C_{m-1}$ and $D$ are the $(m-1)$-th and $m$-th partition set. Now, in order to obtain a partial ${\rm b}^{\ast}$-coloring, assign color $j$ to a vertex if it belongs to the $j$-th partition set. The situation is illustrated in Figure \ref{fig2}, where the vertices of $H$ are vertices with circled labels and the most top vertex represents vertex $u$. A color is also assigned to each vertex with circled label (i.e. the vertices of $H$). As Figure \ref{fig2} shows colors $1, \ldots, i-1,i+1, \ldots, m$ appear in the neighborhood of vertex $u_i$, for each $i=1, \ldots, m$. Hence, $u_i$ can be assigned color $i$ in a ${\rm b}^{\ast}$-coloring of $H$. By assigning colors $1, \ldots, m$ to $u_1, \ldots, u_m$, color $m+1$ is assigned to $u$. The result is a ${\rm b}^{\ast}$-coloring of $H\cup \{u, u_1, \ldots, u_m\}$ with $m+1$ colors. This ${\rm b}^{\ast}$-coloring is depicted in Figure \ref{fig2} for the case $m=4$. Since $G$ is ${\rm b}^{\ast}$-monotonic then ${\rm b}^{\ast}(G)=m+1$. This completes the proof.
\end{proof}

%%%%%%%%%%%%%%%%%%%%%%%%%%%%%%%%%%%%%%%%%%%%%%%%%%%%%%%%%%%%%%%%%%%%%%%%%%%%%%%%%%%%%%%%%%%%

\section{More results on z$(G)$}

\noindent Given a graph $G$, a vertex $u\in V(G)$ and an integer $r\geq 0$, define a ball of radius $r$ centered at $u$ as $B(u,r)=\{v\in V(G):d_G(u,v)\leq r\}$, where $d_G(u,v)$ is the distance i.e. the length of smallest path between $u$ and $v$ in $G$. The following result obtained in \cite{Z3} proves that the Grundy number has a locality property.

\begin{prop}(\cite{Z3})
\noindent Let $G$ be a graph, $c$ be a Grundy-coloring of $G$ using $k$ colors. Let $u$ be a vertex in $G$ such that $c(u)=k$. Let also $H$ be a subgraph of $G$ with minimum cardinality such that $u\in H$ and the restriction of $c$ on $V(H)$ is a Grundy-coloring of $H$ with $k$ colors. Then $V(H)\subseteq B(u, k-1) \subseteq B(u, d_G(u))$.\label{Grunlocal}
\end{prop}

\noindent For a vertex $u$ in a graph $G$ define a new notation ${\bf G^+(u)=G[B(u,d_G(u)+1)]}$.

\noindent The following proposition shows that z-chromatic and ${\rm b}^{\ast}$-chromatic numbers satisfy some locality properties.

\begin{prop}

\noindent $(i)$ For any graph $G$, $z(G)\leq \max_{u\in V(G)} z(G^+(u))$.

\noindent $(ii)$ Equality holds in $(i)$ if $G$ is {\rm z}-monotonic.

\noindent $(iii)$ For any graph $G$, ${\rm b}^{\ast}(G)\leq \max_{u\in V(G)} {\rm b}^{\ast}(G[B(u,2)])$.
\label{zlocal}
\end{prop}

\noindent \begin{proof}
Set for simplicity $M={\max}_{u\in V(G)} z(G^+(u))$. We prove $z(G)\leq M$. Let $c$ be a $z$-coloring of $G$ with $k=z(G)$ colors. Let $u$ be a nice vertex of color $k$ in $c$ and $w_1, \ldots, w_{k-1}$ a set of b-vertex neighbors of $u$ such that $c(w_i)=i$, for each $i=1, \ldots, k-1$. We have $k-1 \leq d_G(u)$. Corresponding to each $w_i$ as a b-vertex of color $i$ and each $j\in \{1, \ldots, i-1, i+1, \ldots, k-1\}$, there exists a neighbors $w^j_i$ of $w_i$ such that $w^j_i$ has color $j$ in the Grundy-coloring $c$. Let $H^j_i$ be a smallest induced subgraph of $G$ such that $w^j_i \in H^j_i$ and the restriction of $c$ on $H^j_i$ is a Grundy-coloring of $H^j_i$ with $j$ colors. Applying Proposition \ref{Grunlocal} for $w^j_i$ and $H^j_i$ we obtain that $V(H^j_i) \subseteq B(w^j_i, j-1)$. Note that $j-1\leq k-2$ and so $V(H^j_i) \subseteq B(w^j_i, k-2)$. It follows that
$$\forall ~ i\in \{1, \ldots, k-1\}~,~ j\in \{1, \ldots, i-1, i+1, \ldots, k-1\}:V(H^j_i) \subseteq B(w_i, k-1) \subseteq B(u, k).$$
\noindent We conclude that $G[B(u,k)]$ and hence $G[B(u,d(u)+1)]$ has a $z$-coloring with $k$ colors. Hence, $z(G)=k \leq z(G^+(u)) \leq M$. This proves $(i)$.

\noindent To prove $(ii)$, note that if $G$ is z-monotonic then a z-coloring of $G^+(u)$ using $z(G^+(u))$ colors is extended to a z-coloring of $G$ with the same number of colors. Proof of $(iii)$ is similar to the proof of $(i)$ and omitted.
\end{proof}

\noindent It was proved in \cite{Z2} that for each $k\geq 1$, there exists a unique tree $R_k$ with a root $u$ such that a rooted tree $T$ with root $w$ admits a z-coloring using $k$ colors in which $w$ receives color $k$ if and only if $T$ contains a subtree isomorphic to $R_k$ such that $u$ is mapped to $w$. It is easily seen by the construction of $R_k$ in \cite{Z2} that the distance of vertices in $R_k$ ($k>3$) from $u$ is at most $k$. In the following we use a result of Varma and Reyner \cite{VR}. Let $T_1$ and $T_2$ be two arbitrary rooted trees with $|V(T_1)|=n$. There exists an $\mathcal{O}(n^{2.5})$ time algorithm which determines whether $T_2$ is isomorphic to some subgraph of $T_2$.

\begin{prop}
Let $G$ be a graph on $n$ vertices and of girth at least $2\Delta(G)+4$. Then $\max_{u\in V(G)} z(G^+(u))$ is determined by an $\mathcal{O}(n\Delta^6)$ time algorithm, where $\Delta=\Delta(G)$.\label{zpoly}
\end{prop}

\noindent \begin{proof}
Let $u$ be an arbitrary vertex in  $G$. Then the girth of $G$ is at least $2d_G(u)+4$. If follows that $G^+(u)$ is an induced tree in $G$. Hence, ${\rm z}(G^+(u))=\max \{k: R_k \unlhd G^+(u)\}$, where $\unlhd$ denotes induced subgraph. We have $|V(G^+(u))|=\mathcal{O}(\Delta^2)$. By the above paragraph $R_k \unlhd G^+(u)$ is decided using $\mathcal{O}(\Delta^5)$ time steps. Since $k\leq \Delta+1$ then ${\rm z}(G^+(u))$ is determined with complexity $\mathcal{O}(\Delta^6)$ and ${\max}_{u\in V(G)} z(G^+(u))$ with time complexity $\mathcal{O}(n\Delta^{6})$.
\end{proof}

%%%%%%%%%%%%%%%%%%%%%%%%%%%%%%%%%%%%%%%%%%%%%%%%%%%%%%%%%%%%%%%%%%%%%%%%%%%%%%%%%%%%%%%%%%%%%%

\noindent To obtain and study integer programming models for the b-coloring are the subject of many papers e.g. \cite{KM}. In the following, we present a 0-1 programming model for z$(G)$. A model for ${\rm b}^{\ast}(G)$ is easily obtained by removing the Grundy-coloring constrain from the formulation for z$(G)$. By Proposition \ref{new-des}, z$(G)$$=$z$'(G)$ which is the maximum number of colors in a Grundy-coloring $G$ which contains a nice vertex. The following formulation is a programming model for z$'(G)$. Let $G$ be a given graph and set a color set $C=\{1, \ldots, m^{\ast}(G)+1\}$.

\begin{enumerate}
  \item $\max \mathcal{Z}=\sum_{c\in C} m_c$
  \item ~~~~~~$\sum_{c\in C} x_{vc} =1$,~~~~~~~~ $\forall v\in V$
  \item ~~~~~~$x_{vc}+x_{uc} \leq m_c$,~~~  $\forall c\in C$, $\forall u,v: uv\in E$
  \item ~~~~~~$x_{vc}\leq m_c$,~~~~~~~~~~~ $\forall v\in V$, $\forall c\in C$
  \item ~~~~~~$m_c \leq \sum_{v\in V} x_{vc}$, ~~~~~~$\forall c\in C$
  \item ~~~~~~$x_{vc'} \leq \sum_{u\in N(v)} x_{uc}$, ~~$\forall v\in V$,  $\forall c,c' \in C$: $c<c'$
  \item ~~~~~~$m_{c'} \leq m_c$, ~~~~~~~~~~~$\forall c, c' \in C$: $c<c'$
  \item ~~~~~~$z_{vc}\leq \sum_{u\in N(v)} x_{ud}$, ~~~~$\forall v\in V$, $\forall c,d \in C$: $c\not= d$
  \item ~~~~~~$z_{vc} \leq x_{vc}$, ~~~~~~~~~~~ $\forall v\in V$, $\forall c\in C$
  \item ~~~~~~$\zeta_{vc}\leq \sum_{u\in N(v)} z_{ud}$, ~~~~~$\forall v\in V$, $\forall c,d \in C$: $c\not= d$
  \item ~~~~~~$\zeta_{vc} \leq z_{vc}$, ~~~~~~~~~~~~ $\forall v\in V$, $\forall c\in C$
  \item ~~~~~~$1 \leq \sum_{v\in V, c\in C} \zeta_{vc}$
  \item ~~~~~~$x_{vc},~ m_c,~ z_{vc},~ \zeta_{vc} \in \{0,1\}$
\end{enumerate}

\begin{prop}
The objective function $\mathcal{Z}$ in the programming model determines ${\rm z}(G)$ correctly.\label{prog}
\end{prop}

\noindent \begin{proof}
\noindent For each $c\in C$, $m_c=1$ if and only if $c$ is assigned for at least one vertex in $G$. For $c\in C$ and $v\in V$, $x_{vc}=1$ if and only if $c$ is assigned as a color to vertex $v$. Constraint (2) guarantees that each vertex is received exactly one color. Constraint (3) guarantees that the assignment is proper coloring of $G$. Hence, $\{x_{vc}: x_{vc}\not= 0\}$ defines a proper vertex coloring of $G$ using $\max \{c: m_c\not= 0\}$ colors. We prove that this coloring is ${\rm z}'$-coloring. Constraint (5) asserts that if $m_c=1$ for some color $c$, then there exists a vertex having color $c$. Constraint (6) guarantees that the coloring has Grundy property. Constraint (7) means that $\{c:m_c\not=0\}$ is a continuous set of integers. This is a necessary condition since the coloring has Grundy property and also the task is to enumerate the number of distinct colors. Constraint (8) guarantees that if $z_{vc}=1$ then $v$ is b-vertex of color $c$. Constraint (9) regulates the values $z_{vc}$ in terms of $x_{vc}$. Constraint (10) guarantees that if $\zeta_{vc}=1$ then for each color $c'\not= c$, $v$ has a b-vertex neighbor of color $c'$, i.e. $v$ is a nice vertex of color $c$. Constraint (11) means that a nice vertex of color $c$ is surely b-vertex of color $c$. Constrain (12) guarantees that there exists a color $c$ and a vertex $v$ such that $v$ is nice vertex of color $c$. The objective function $\mathcal{Z}$ in (1) maximizes the number of used colors.
Note that every graph $G$ admits a coloring with $\chi(G)$ colors such that contains a nice vertex. So the space of feasible solutions for the program is non-empty. It follows that the program finds the maximum number of colors $k$ such that there exists a proper Grundy-coloring with $k$ colors containing a nice vertex. Hence, $k={\rm z}'(G)$ and by Proposition \ref{new-des}, $k={\rm z}(G)$.
\end{proof}

\begin{cor}
If we remove the constraint (6) from the model in Proposition \ref{prog} then the result is a 0-1 programming model for ${\rm b}^{\ast}(G)$.
\end{cor}

\section{Suggestions for further researches}

\noindent As z-coloring and ${\rm b}^{\ast}$-coloring are new concepts, there are many unexplored problems and research areas involving these types of colorings. A main goal is to prove the algorithmic or analytic advantages of these colorings in comparison to the First-Fit and b-colorings. We specify the following problems.

\noindent {\bf Problem 1.} Study of some chromatic invariants in terms of the so called atom graphs was initiated for Grundy chromatic number in \cite{Z1} and then for b-chromatic number in \cite{EGT} and recently for z-chromatic number in \cite{Z2}. Define b$^{\ast}$-atoms for b$^{\ast}$-colorings and provide a constructive method to generate them. As b$^{\ast}$-atoms should be connected their structure will surely be simpler than the b-atoms.

\noindent {\bf Problem 2.} Proposition \ref{bsparse} asserts that ${\rm b}^{\ast}(G) =\omega(G)$ for $P_4$-free graphs. Does there exist a function $f(.)$ such that for any $P_5$-free graph $G$, ${\rm b}^{\ast}(G)\leq f(\omega(G))$? Kierstead et al. in \cite{KPT} proved that there exists an exponential function $g(.)$ such that $\Gamma(G)\leq g(\omega(G))$ for any $P_5$-free graph $G$. It is interesting to explore a function $f(.)$ of a lower magnitude than $g(.)$ such that ${\rm z}(G)\leq f(\omega(G))$, for $P_5$-free graphs $G$. It was proved in \cite{Z2} using z-atoms that ${\rm z}(G)\leq 3$, for $(P_5,K_3)$-free graphs. What about the z-chromatic and ${\rm b}^{\ast}$-chromatic numbers of $(P_5,K_4)$-free graphs?

\noindent {\bf Problem 3.} ${\rm b}^{\ast}(G)$ is polynomially solvable for block graphs, cacti and $P_4$-sparse graphs. Find more classes of graphs for which the problem is polynomially solvable.

\noindent {\bf Problem 4.} What is the smallest $g$ such that every graph of girth at least $g$ is ${\rm b}^{\ast}$-monotonic? Similar question for z-monotonicity is interesting.

\noindent {\bf Problem 5.} Coloring heuristics are usually compared using their outputs on random graphs or on DIMACS benchmarks. Book \cite{L} reports the results for Greedy, DSATUR, Recursive Largest First (RLF) and some other heuristics. Iterated Greedy (IG) defined by J. Culberson (see \cite{CL}) has been compared with heuristics such as TabuCol (\cite{HW}) in \cite{CL}. As we mentioned earlier, a heuristic derived from z-coloring is defined in \cite{Z2} and denoted by IZ. An interesting research is to compare the performance of IZ with heuristics such as IG, TabuCol, DSATUR and RLF. Theoretically, IZ is more optimal than IG but no experimental comparison has been performed.

%%%%%%%%%%%%%%%%%%%%%%%%%%%%%%%%%%%%%%%%%%%%%%%%%%%%%%%%%%%%%%%%%%%%%%%%%%%%%%%%%%%%%%%%%%%%%%%%%%%%%%%%%%%%

\end{document}